\documentclass[11pt]{article}
\usepackage{amsmath}
\usepackage{amsthm}
\usepackage{latexsym}
\usepackage[noblocks]{authblk}
\usepackage{amssymb,amsmath}  
\usepackage{multimedia}
\usepackage{subfigure}
\usepackage{enumerate}
\usepackage{cite}
\usepackage{authblk} 
\usepackage{mathrsfs} 
\usepackage{epsfig}
\usepackage{mathdots}
\usepackage{caption}
\usepackage{color}
\usepackage{pgfpages,tikz,tikz-cd,pgfkeys,pgfplots}
\usetikzlibrary{arrows,positioning,matrix,fit,backgrounds,shapes}

\usetikzlibrary{fit,matrix}
\tikzset{
mN/.style = {
    draw=#1, semithick, inner sep=0pt}
             }

\usetikzlibrary{arrows,positioning,matrix,fit,backgrounds,shapes,shapes.geometric, intersections}
\usetikzlibrary{shapes.callouts,decorations.pathmorphing}
\usetikzlibrary{decorations.markings} 
\tikzstyle{vertex}=[circle, draw, inner sep=0pt, minimum size=6pt] 

\definecolor{LemonChiffon}{rgb}{100, 98, 80}
\definecolor{myblue}{rgb}{0,0.4,0.8}
\definecolor{orange}{rgb}{1, 0.4, 0}
\definecolor{mygreen}{rgb}{0, 0.8, 0.2}
\definecolor{myred}{rgb}{204, 0, 0}
\definecolor{violet}{RGB}{0.4,0.2,1}
\definecolor{brown}{rgb}{0.6, 0.4, 0}

\newtheorem{theorem}{Theorem}[section]
\newtheorem{lemma}[theorem]{Lemma}

\newtheorem{proposition}[theorem]{Proposition}
\newtheorem{corollary}[theorem]{Corollary}

\theoremstyle{definition}

\begin{document}

\title{Competition graphs of degree bounded digraphs}
\author[1]{Hojin Chu}
\author[1]{Suh-Ryung Kim}

\affil[1]{\footnotesize Department of Mathematics Education, Seoul National University, Seoul 08826, Rep. of Korea}
\affil[ ]{\footnotesize\textit{ghwls8775@snu.ac.kr, srkim@snu.ac.kr}}
\date{}

\maketitle
\begin{abstract}
If each vertex of an acyclic digraph has indegree at most $i$ and outdegree at most $j$, then it is called an $(i,j)$ digraph, which was introduced by Hefner~{\it et al.}~(1991).  
Whereas Hefner~{\it et al.} characterized $(i,j)$ digraphs whose competition graphs are interval, characterizing the competition graphs of $(i,j)$ digraphs is not an easy task.
In this paper, we introduce the concept of $\langle i,j \rangle$ digraphs, which relax the acyclicity condition of $(i,j)$ digraphs, and study their competition graphs.
By doing so, we obtain quite meaningful results.
Firstly, we give a necessary and sufficient condition for a loopless graph being an $\langle i,j \rangle$ competition graph for some positive integers $i$ and $j$.
Then we study on an $\langle i,j \rangle$ competition graph being chordal and present a forbidden subdigraph characterization.
Finally, we study the family of $\langle i,j \rangle$ competition graphs, denoted by $\mathcal{G}_{\langle i,j \rangle}$, and identify the set containment relation on $\{\mathcal{G}_{\langle i,j \rangle}\colon\, i,j \ge 1\}$.
\end{abstract}



\section{Introduction}
In this paper, all the digraphs are assumed to have neither loops nor parallel arcs. 
For all undefined graph-theoretical terminology, we follow~\cite{Bondy2010graph}.

In 1968, Cohen~\cite{Cohen1968} introduced the notion of competition graphs in connection with competition in the food web of ecosystems.
Given a digraph $D$,  the \textit{competition graph} of $D$, denoted by $C(D)$, has the same vertex set as $D$ and an edge $uv$ if and only if $u$ and $v$ have a common out-neighbor in $D$.
A competition graph of a digraph is simple.
In this vein, all the graphs are assumed to be simple throughout the paper.
The competition graph has various applications in coding, radio transmission, and modeling of complex economic systems and its variants were introduced and have been studied (see \cite{kim1993competition}, \cite{JRIchard1989food},\cite{Raychaudhuri1985},\cite{cho2000m},\cite{Choi2023},\cite{lee2017phylogeny}, \cite{Eoh2021chordalphylogeny}).
Especially, we note that \cite{lee2017phylogeny} and \cite{Eoh2021chordalphylogeny} are research results derived from recent studies on phylogeny graph, which is a variant of competition graph, of a degree bounded digraph.

When Cohen introduced the notion of competition graphs, he observed that most food web is acyclic, that is, there is no directed cycle in the digraph modeled by the food web.
Thus mostly competition graphs of acyclic digraphs were studied in the early stages of research.
In 1982, Opsut~\cite{Opsut1982} showed that the problem of determining whether or not an arbitrary graph is the competition graph of some acyclic digraph is NP-complete. 
Since it is not easy to characterize the competition graphs of acyclic digraphs in general, researchers have been attempted to study the competition graphs of restricted acyclic digraphs or arbitrary digraphs.
For example, 
Dutton and Brigham~\cite{RonaldD1983} characterized the competition graphs of arbitrary digraphs.
Roberts and Steif~\cite{Roberts1983} characterized the competition graphs of loopless digraphs. 
Lundgren and Maybee~\cite{Richard1983} characterized the graphs which can be made into the competition graphs of acyclic digraphs by adding $m$ isolated vertices but it is not possible to be the competition graphs of acyclic digraphs if less than $m$ isolated vertices are added.
While Hefner~{\it et al.}~\cite{Hefner1991ij} characterized certain degree bounded acyclic digraphs whose competition graphs are interval, it is not easy to characterize even the competition graphs of degree bounded acyclic digraphs. 
In this paper, we characterize the competition graphs of degree bounded loopless digraphs.

In 1991, Hefner~{\it et al.} introduced the notion of $(i,j)$ digraphs and studied the competition graphs of $(i,j)$ digraphs.
The \textit{$(i,j)$ digraph} is an acyclic digraph such that $d^{-}(x) \le i$ and $d^{+} (x) \le j$ for every vertex $x$, where $d^{-}(x)$ and $d^{+}(x)$ denote the indegree and the outdegree of a vertex $x$, respectively.
We expand the notion of $(i,j)$ digraphs by allowing some directed cycles and call it an {\it $\langle i,j \rangle$ digraph}.
The \textit{$\langle i,j \rangle$ digraph} is a loopless digraph satisfying $d^{-}(x) \le i$ and $d^{+} (x) \le j$ for every vertex $x$.
By definition, an $\langle i,j \rangle$ digraph might not be acyclic.
We say that a graph is an {\it $\langle i,j \rangle$ competition graph} if it is the competition graph of some $\langle i,j \rangle$ digraph.
By definition, note that an $\langle i,j \rangle$ competition graph is a $\langle k, l \rangle$ competition graph for each integers $i$, $j$, $k$, and $l$ satisfying $i\le k$ and $j\le l$.

In Section~\ref{sec:2}, we give necessary and sufficient conditions for a loopless graph being an $\langle i,j \rangle$ competition graph for some positive integers $i$ and $j$ (Theorems~\ref{thm:edge clique cover} and \ref{thm:edge clique cover_2}).
We also deduce an interesting fact that the existence of certain $\langle i,j \rangle$ competition graph is equivalent to the existence of a balanced incomplete block design for some variables (Proposition~\ref{prop:clique size}).
In addition, we give a sufficient condition for an $\langle i,j \rangle$ competition graph being chordal in the aspect of subdigraph restriction (Proposition~\ref{prop:chordal if}).
Especially, it becomes a necessary and sufficient condition when $j=2$ (Theorem~\ref{thm:i,2 chordal}).
In Section~\ref{sec:3}, we denote the family of $\langle i,j \rangle$ competition graphs by $\mathcal{G}_{\langle i,j \rangle}$ and take a look at set containment between $\mathcal{G}_{\langle i,j \rangle}$ and $\mathcal{G}_{\langle k,l \rangle}$ for given positive integers $i$, $j$, $k$, and $l$ with $(i,j) \neq (k,l)$ (Theorem~\ref{thm:containment}).
To show $\mathcal{G}_{\langle i,j \rangle}\not \subseteq \mathcal{G}_{\langle k,l \rangle}$ for specific $i$, $j$, $k$, and $l$, we use quite novel idea to construct a graph belonging to $\mathcal{G}_{\langle i,j \rangle}$ but not belonging to $\mathcal{G}_{\langle k,l \rangle}$ (Propositon~\ref{prop:tool3}).

\section{$\langle i, j \rangle$ competition graphs}\label{sec:2}

Given a graph $G$, a {\it clique} $C$ is a subgraph of $G$ in which any two vertices are adjacent.
If there is no possibility of confusion, we will use the same symbol ``C" for the set of vertices of C.
An {\it edge clique cover} of a graph $G$ is the collection of cliques in $G$ that covers all the edges in $G$.

Let $A$ be a finite set and let $\mathcal{A}=\{A_1, \ldots, A_m\}$ be a collection of subsets of $A$.
A {\it system of distinct representatives} (SDR) of $\mathcal{A}$ is a collection of distinct elements $a_1, \ldots, a_m$ such that $a_i \in A_i$ for each $1\le i \le m$.

Roberts and Steif gave a 
necessary and sufficient condition for a graph being the competition graph of a loopless digraph as follows.

\begin{theorem}[Roberts and Steif\cite{Roberts1983}]
A graph $G$ is the competition graph of a digraph which has no loop if and only if there is an edge clique cover $\{C_1, C_2, \ldots, C_p\}$ of $G$ such that $\{V(G)-C_1,V(G)-C_2, \ldots, V(G)-C_p\}$ has a system of distinct representatives.
\end{theorem}

An $\langle i,j \rangle$ competition graph obviously satisfies the condition given in the above theorem.
As a matter of fact, it satisfies three more conditions stated in the following theorem and those four conditions in the theorem are actually sufficient to guarantee that a graph is an $\langle i,j \rangle$ competition graph.

\begin{theorem}\label{thm:edge clique cover}
A graph $G$ is an $\langle i,j \rangle$ competition graph for some positive integers $i$ and $j$ if and only if there is an edge clique cover $\mathcal{C}=\{C_1, \ldots, C_p\}$ of $G$ satisfying the following conditions:
\begin{enumerate}[(i)]
	\item $|C_t| \le i$ for each $1\le t \le p$;
	\item each vertex of $G$ belongs to at most $j$ cliques in $\mathcal{C}$;
	\item $\{V(G)-C_1, \ldots, V(G)-C_p\}$ has an SDR;
	\item $p\le |V(G)|$.
\end{enumerate}
\end{theorem}
\begin{proof}
	To show the ``only if" part, suppose that $G$ is the competition graph of an $\langle i,j \rangle$ digraph $D$.
	Let $\mathcal{C}=\{ N^-_D(v)\colon\, v\in V(D) \text{ with } N^-_D(v)\neq \emptyset \}=\{N^-_D(v_1), \ldots, N^-_D(v_p)\colon\, v_1,\ldots, v_p \in V(D)\}$.
	By the definition of competition graph, $\mathcal{C}$ is an edge clique cover of $G$.
	By the fact that $G$ is an $\langle i,j \rangle$ competition graph, $\mathcal{C}$ satisfies the conditions (i) and (ii).
	Moreover, since $D$ is loopless, $\{v_t \colon\, 1\le t \le p\}$ is an SDR of $\{V(G)-N^-_D(v_t)\colon\, 1\le t \le p \}$ and so 
	$\mathcal{C}$ satisfies the condition (iii).	
	Then, by the condition (iii), there are $p$ distinct vertices of $G$, which happen to be distinct representatives.
	Thus $p \le |V(G)|$ and so the condition (iv) holds.
	
	To show the ``if" part, suppose that there is an edge clique cover $\mathcal{C}=\{C_1, \ldots, C_p\}$ of given graph $G$ satisfying the conditions (i), (ii),  (iii), and (iv).
	Let $\{v_1, \ldots, v_p\}$ be an SDR of $\{V(G)-C_1, \ldots, V(G)-C_p\}$.
	Now, we define a digraph $D$ with the vertex set $V(G)$ and the arc set
	\[\bigcup_{t=1}^p \{(v,v_t) \colon\, v\in C_t\}.\]
	Since $v_t \notin C_t$ for each $1\le t \le p$, $D$ has no loop.
	Take a vertex $v$ in $D$.
	Then $d_D^+(v)\le j$ by the condition (ii).
	If $v \neq v_t$ for any $1\le t \le p$, then $N_D^-(v)=\emptyset$.
	If $v = v_t$ for some $1\le t \le p$, then $N_D^-(v)=C_t$ and so $d_D^-(v) \le i$ by the condition (i).
	Thus $D$ is an $\langle i,j \rangle$ digraph.
	Moreover, $u$ and $v$ are adjacent in $G$ if and only if $\{u,v\} \subseteq C_t$ for some $1\le t \le p$ if and only if $(u,v_t)$ and $(v,v_t)$ are arcs of $D$ for some $1\le t \le p$ if and only if $u$ and $v$ are adjacent in the competition graph of $D$.
	Therefore the competition graph of $D$ is $G$.
\end{proof}

\begin{corollary}\label{cor:edge clique cover}
Let $G$ be an $\langle i,j \rangle$ competition graph for some positive integers $i$ and $j$.
Then the following are true:
\begin{center}
	(a) $|E(G)|\le \frac{(i-1)i}{2}|V(G)|$; \quad (b) $G$ is  $K_{1,j+1}$-free; \quad (c) $\Delta(G)\le j(i-1)$,
\end{center}
where $\Delta(G)$ denotes the maximum degree of $G$, that is, the largest vertex degree of $G$.
\end{corollary}
\begin{proof}
	By Theorem~\ref{thm:edge clique cover}, there is an edge clique cover $\mathcal{C}:=\{C_1, C_2, \ldots, C_p\}$ of $G$ which satisfies the conditions (i), (ii), (iii), and (iv) of the theorem.
	By the condition (iv), $p \le |V(G)|$.
	Thus
\[		|E(G)|= \left|\bigcup_{t=1}^{p}E(C_t)\right| \le \sum_{t=1}^{p}{|C_t| \choose 2} \le \sum_{t=1}^{p}{i \choose 2} \le \sum_{t=1}^{|V(G)|}{i \choose 2}= \frac{(i-1)i}{2}|V(G)|
\]
where the second inequality holds by the condition (i).
Therefore (a) holds.

To show (b) by contradiction, suppose that $G$ has an induced subgraph $H$ isomorphic to $K_{1,j+1}$.
	We denote the center of $H$ by $v$.
	Then the neighbors of $v$ in $H$ belong to distinct cliques in $\mathcal{C}$ and so $v$ belongs to at least $j+1$ distinct cliques, which contradicts the condition (ii).
	
To show (c), take a vertex $v$.
Let $C_{p_1}, \ldots, C_{p_s}$ be the cliques in $\mathcal{C}$ containing $v$.
Then $s \le j$ by the condition (ii).
Now
\[d_G(v)=|N_G(v)|=\left|\bigcup_{t=1}^s (C_{p_t}-\{v\}) \right|\le \sum_{t=1}^s \left(|C_{p_t}|-1 \right) \le \sum_{t=1}^s \left(i-1 \right) \le j(i-1) \]
where the second inequality holds by the condition (i).
Therefore (c) holds.
\end{proof}

\begin{lemma}\label{lem:complete graph}
For each integer $n\ge 3$, 
a complete graph $K_n$ is an $\langle i,j \rangle$ competition graph if $i\ge n-1$ and $j\ge 2$.
Further, 
$K_n$ is not a $\langle k, 1\rangle$ competition graph for any integer $k\ge 1$.
\end{lemma}
\begin{proof}
	Fix an integer $n \ge 3$.
	We define a digraph $D$ with the vertex set $\{v_1,\ldots,v_n\}$ and the arc set
	\[\{(v_t, v_n) \colon\, 1\le t \le n-1\} \cup \{(v_t, v_1) \colon\, 2 \le t \le n\} \cup \{(v_1,v_2), (v_n, v_2)\}.\]
	Then it is easy to see that $D$ is an $\langle n-1,2 \rangle$ digraph and its competition graph is isomorphic to $K_n$.
	Therefore a complete graph $K_n$ is an $\langle n-1,2 \rangle$ competition graph and so it is an $\langle i,j \rangle$ competition graph for any integers $i\ge n-1$ and $j\ge 2$.	
	
	Let $G$ be the competition graph of a $\langle k,1 \rangle$ digraph $D$ of order $n$ for an integer $k\ge 1$.
	Take two adjacent vertices $u$ and $v$ in $G$.
	Then there is a common out-neighbor $w$ of $u$ and $v$.
	Since $D$ is a $\langle k,1 \rangle$ digraph without loops, the vertices $u$, $v$, and $w$ are distinct and \[N_D^+(u)=N_D^+(v)=\{w\}.\]
	Thus $u$ and $w$ have no common out-neighbor and hence are not adjacent in $G$.
	Therefore a complete graph $K_n$ is not a $\langle k, 1\rangle$ competition graph for any integer $k\ge 1$.
\end{proof}

\begin{lemma}[Roberts and Steif\cite{Roberts1983}]\label{lem:subsets}
Let $k$ be a positive integer and $S_1, \ldots, S_{k+1}$ be subsets of $\{1, \ldots, k\}$.
Then there is an integer $t \in \{1,2, \ldots, k+1\}$ such that $S_t\subseteq \bigcup_{i\neq t} S_i$.
\end{lemma}

In Theorem 2.2, the conditions (i), (ii), (iii), and (iv) are given for a graph to be an $\langle i, j\rangle$ competition graph. 
However, when conditions (i), (ii), and (iv) hold, the condition (iii) is equivalent to the following statement:
$G\not\cong K_2$ and, if $j=1$, then $G$ is not a nontrivial complete graph.
In this context, for a graph $G$ and positive integers $i$ and $j$, we denote by $\mathbf{C}(G,i,j)$ the collection of edge clique cover $\mathcal{C}=\{C_1, \ldots, C_p\}$ of $G$ satisfying the conditions (i), (ii), and (iv) of Theorem~\ref{thm:edge clique cover}.
Then we obtain another necessary and sufficient condition for a graph to be an $\langle i, j\rangle$ competition graph as follows. 
To this end, we need Hall's Marriage Theorem.
\begin{theorem}[Hall\cite{hall1935}]\label{thm:MC}
	The family $\mathcal{A}=(A_1,A_2,\ldots,A_n)$ of sets has a system of distinct representatives
	if and only if 
	\begin{center}
	\begin{minipage}{0.85\textwidth}
	(Hall's marriage condition)
for each $1\leq k\leq n$ and each choice of $k$ distinct indices $i_1,i_2,\ldots,i_k$ from $\{1,2,\ldots,n\}$,
\[\left|A_{i_1} \cup A_{i_2} \cup \cdots \cup A_{i_k}\right| \geq k.\]
	\end{minipage}
	\end{center}
\end{theorem}

\begin{theorem}\label{thm:edge clique cover_2}
	A graph $G$ is an $\langle i,j \rangle$ competition graph for some positive integers $i$ and $j$ if and only if (i) $G\not\cong K_2$, (ii) $G$ is not a nontrivial complete graph for $j=1$, and (iii) $\mathbf{C}(G,i,j)\neq \emptyset$.
\end{theorem}
\begin{proof}	
	It is easy to see that a complete graph $K_2$ cannot be the competition graph of any digraph unless a loop is allowed.
	Then the ``only if" part is true by Theorem~\ref{thm:edge clique cover} and the ``further" part of Lemma~\ref{lem:complete graph}.

	To show the ``if" part, let $G$ be a graph of order $n$ satisfying the conditions (i), (ii), and (iii).
	Let $\mathcal{C}=\{C_1, \ldots, C_p\}$ be an edge clique cover in $\mathbf{C}(G,i,j)$ with the smallest $w(\mathcal{C}):=\sum_{t=1}^{p} |C_t|$.
	By (iii), $\mathcal{C}$ is well-defined.
	If $\sum_{t=1}^{p} |C_t|=1$, then $G$ is edgeless and so it is the competition graph of a digraph $D=(V,\emptyset)$ which is an $\langle i,j \rangle$ digraph for any positive integers $i$ and $j$.
	
	Suppose $\sum_{t=1}^{p} |C_t| > 1$.
	Then $G$ has an edge.
	If $p=1$ and $G$ has an isolated vertex $v$, then $G$ is the competition graph of an $\langle i,j \rangle$ digraph $D$ whose vertex set is $V(G)$ and arc set is
	\[\{(u,v) \colon\, u \in C_1\}.\] 
	Consider the case where $p=1$ and $G$ has no isolated vertex.
	Then $G$ is a complete graph.
	Thus $j\ge 2$ and $n\ge 3$ by (i) and (ii). 
	Therefore $G$ is an $\langle i,j \rangle$ competition graph by Lemma~\ref{lem:complete graph}.
	Thus we may assume $p\ge 2$.	
	Then $G$ is not complete since $w(\mathcal{C})$ is the smallest and no clique in $\mathcal{C}$ is contained in another clique by the choice of $\mathcal{C}$.
	We will show that $\mathcal{B}:=\{B_1, \ldots, B_p\}$ has an SDR where $B_t=V(G)-C_t$ for each $1\le t \le p$.
	Then we are done by Theorem~\ref{thm:edge clique cover}.
	
	To the contrary, suppose that $\mathcal{B}:=\{B_1, \ldots, B_p\}$ has no SDR.
	Then, by Theorem~\ref{thm:MC}, 
	there is an nonempty subset $\{t_1, \ldots, t_s\}$ of $\{1,2, \ldots, p\}$ satisfying 
\[B_{t_1} \cup \cdots \cup B_{t_s}=\{v_1, \ldots, v_l\}\]
where $v_1, \ldots, v_l$ are vertices of $G$ for some integer $l<s$.
Since $G$ is not complete, $V(G)-C_t \neq \emptyset$ for each $1\le t\le p$.
Thus $l \ge 1$ and so $s \ge 2$.
For each $1 \le a \le s$,
\[C_{t_a} \supseteq C_{t_1} \cap \cdots \cap C_{t_s}=V(G)-(B_{t_1} \cup \cdots \cup B_{t_s})=V(G)-\{v_1, \ldots, v_l\}. \]
Thus, for each $1\le a \le s$, $C_{t_a}=(V(G)-\{v_1, \ldots, v_l\}) \cup Q_{t_a}$ for some $Q_{t_a} \subseteq \{v_1, \ldots, v_l\}$.
Then $Q_{t_a}$ is a subset of $C_{t_a}$ which is a clique, so $Q_{t_a}$ is also a clique for each $1\le a \le s$.
Since $Q_{t_1}, \ldots, Q_{t_s}$ are subsets of $\{v_1, \ldots, v_l\}$ and $s>l$, there is an integer $r \in \{1,2, \ldots, s\}$ such that 
\[Q_{t_r} \subseteq  \bigcup_{i\neq r}Q_{t_i}\]
by Lemma~\ref{lem:subsets}.
Then $\mathcal{C}':= (\mathcal{C}\setminus\{C_{t_r}\}) \cup \{Q_{t_r}\}$ belongs to $\mathbf{C}(G,i,j)$.
To see why, we first observe that $\mathcal{C}'$ consists of cliques of $G$.
For an edge $uv$ in $C_{t_r}$, (1) $\{u,v\} \subseteq V(G)-\{v_1, \ldots, v_l\}$ or (2) $u \in V(G)-\{v_1, \ldots, v_l\}$ and $v \in Q_{t_r}$ or (3) $u \in Q_{t_r}$ and $v \in V(G)-\{v_1, \ldots, v_l\}$ (4) $\{u,v\} \subseteq Q_{t_r}$.
The cases (4) and (1) are taken care of since $Q_{t_r} \in \mathcal{C}'$ and $V(G)-\{v_1, \ldots, v_l\} \subseteq C_{t_a}$ for $1\le a \le s$, respectively.
Consider the cases (2) and (3).
Without loss of generality, we may assume that (2) happens.
Since $Q_{t_r} \subseteq (Q_{t_1}\cup \cdots \cup Q_{t_s})-Q_{t_r}$, there is an integer $a \in \{1,2, \ldots, s\}$ with $a \neq r$ such that $v \in Q_{t_a}$.
Then $v \in C_{t_a}$.
Since $V(G)-\{v_1, \ldots, v_l\} \subseteq C_{t_a}$, $u \in C_{t_a}$.
Thus the edge $uv$ is covered by $C_{t_a} \in \mathcal{C}'$.
Therefore $\mathcal{C}'$ is an edge clique cover of $G$.
Since $Q_{t_r}\subsetneq C_{t_r}$, $\mathcal{C}' \in \mathbf{C}(G,i,j)$ and $w(\mathcal{C}') < w(\mathcal{C})$, which contradicts the choice of $\mathcal{C}$. 
Hence $G$ is an $\langle i,j \rangle$ competition graph.
\end{proof}

We recall some definitions about balanced incomplete block designs.
A {\it design} is a pair $(V, \mathcal{B})$ where $V$ is a finite set of varieties and $\mathcal{B}$ is a collection of subsets of $V$, called {\it blocks} of the design.
A {\it balanced incomplete block design (BIBD)} is a design consisting of $v$ varieties and $b$ blocks such that each block consists of exactly the same number $k$ of varieties where $k<v$; each variety appears in exactly the same number $r$ of blocks; each pair of varieties appears simultaneously in exactly the same number $\lambda$ of blocks.
Such a design is called a {\it $(b,v,r,k,\lambda)$-BIBD}.

The {\it clique number} of a graph $G$, denoted by $\omega(G)$, is the number of vertices in a maximum clique of $G$.

\begin{proposition}\label{prop:clique size}
Let $G$ be an $\langle i,j \rangle$ competition graph for some positive integers $i$ and $j$.
	Then $\omega(G) \le ij-j+1$.
	Further, the upper bound is sharp if and only if there is a $(j(ij-j+1)/i, ij-j+1, j, i, 1)$-BIBD.
\end{proposition}
\begin{proof}	
By Corollary~\ref{cor:edge clique cover}(c), it directly follows that $\omega(G) \le j(i-1)+1$.
	
	There exists an $\langle i,j \rangle$ digraph $D$ whose competition graph is $G$.
	To show the ``only if" part of the ``further" part, suppose that $G$ has a clique $K$ of size $ij-j+1$.
	For any $u\in K$,
	\begin{align*}
		ij-j=d_{G[K]}(u)=\left|\bigcup_{v \in N^+_D(u)}(N^-_D(v)-\{u\}) \cap K\right| \le \sum_{v \in N^+_D(u)}|(N^-_D(v)-\{u\}) \cap K| \\ \le \sum_{v \in N^+_D(u)}|(N^-_D(v)-\{u\})| \le\sum_{v \in N^+_D(u)}(i-1) \le (i-1)j
	\end{align*}
	where $G[K]$ denotes the subgraph of $G$ induced by $K$.
	Thus the inequalities above become equalities.
	The first inequality being equality implies that	for any distinct $u$ and $v$ in $K$, 
	\begin{equation}\label{eq:1}
		\left|N_D^+(u)\cap N_D^+(v)\right|=1.
	\end{equation}
	It follows from the conversions of the second and the third inequalities that for any vertex $v \in \bigcup_{u \in K}N^+_D(u)$,
	\begin{equation}\label{eq:2}
		N^-_D(v) \subseteq K \quad \text{ and } \quad d_D^-(v)=i,
	\end{equation}
	respectively.
	By the conversion of the last inequality, it is immediately true that for each vertex $u$ in $K$,
	\begin{equation}\label{eq:3}
		d_D^+(u)=j.
	\end{equation}
	Set \[L=\{N^-_D(v)\colon\, v \in N^+_D(u)\text{ for some } u \in K\}.\]
	We regard the vertices in $K$ (resp.\ the elements in $L$) as the varieties (resp.\ the blocks).
	Then, by \eqref{eq:1}, \eqref{eq:2}, and \eqref{eq:3},  we have a $(|L|,ij-j+1,j,i,1)$-BIBD.
	It is well-known that for a $(b,v,r,k,\lambda)$-BIBD, $bk=vr$.
	Therefore $|L|=j(ij-j+1)/i$.
	Thus there is a $(j(ij-j+1)/i, ij-j+1, j, i, 1)$-BIBD.
	
	To show the ``if" part, suppose that there is a $(j(ij-j+1)/i, ij-j+1, j, i, 1)$-BIBD.
	Let $\{x_1,x_2,\ldots, x_{ij-j+1}\}$ be the set of varieties and $B_1, \ldots, B_{j(ij-j+1)/i}$ be the blocks of the design.
	Now we define a digraph $D$ such that the vertex set is
	\[\{u_1,u_2, \ldots, u_{ij-j+1}\}\cup \{v_1,v_2, \ldots, v_{j(ij-j+1)/i}\}\] 
	and $D$ has an arc $(u_k,v_l)$ if and only if $x_k \in B_l$ for some integers $1\le k \le ij-j+1$ and $1\le l \le j(ij-j+1)/i$.
	Then, since each variety appears $j$ blocks and each block has $i$ varieties, $D$ is an $\langle i,j\rangle$ digraph.
	Moreover, since each pair of varieties is contained in a unique block, $\{u_1,u_2,\ldots, u_{ij-j+1}\}$ forms a clique in the competition graph of $D$.
	Therefore $\omega(C(D)) \ge ij-j+1$.
	As we have shown $\omega(C(D)) \le ij-j+1$, $\omega(C(D))=ij-j+1$.
	\end{proof}

\begin{corollary}
The following are true for any positive integers $i$ and $j$.
	\begin{enumerate}[(a)]
		\item $\omega(G) < ij-j+1$ for an $\langle i,j \rangle$ competition graph $G$ if  $i>j$.
		\item There is a $\langle 2,j \rangle$ competition graph with clique number $j+1$.
		\item There is a $\langle 3,j \rangle$ competition graph with clique number $2j+1$ if and only if $j=1$ or $j=3n$ or $j=3n+1$ for some integer $n\ge 1$.
	\end{enumerate}
\end{corollary}
\begin{proof}
	Let $G$ be an $\langle i,j \rangle$ competition graph with $i>j$.
	By Proposition~\ref{prop:clique size}, $\omega(G) \le ij-j+1$.
	To the contrary, suppose $\omega(G) = ij-j+1$.
	Then, by the ``further" part of Proposition~\ref{prop:clique size}, there is a $(j(ij-j+1)/i, ij-j+1, j, i, 1)$-BIBD.
	Now, by Fisher's inequality, $j \ge i$, which is a contradiction.
	Therefore $\omega(G) < ij-j+1$ and the part (a) is true.
	
	There is a $(j(j+1)/2, j+1, j, 2, 1)$-BIBD for each positive integer $j$ since the two-element subsets of a set of $j+1$ elements correspond blocks.
	Thus, by the ``further" part of Proposition~\ref{prop:clique size}, the part (b) is true.
	
	It is well-known that a Steiner triple system $S(2,3,n)$ exists if and only if $n=3$ or $n=6k+1$ or $n=6k+3$ for some integer $k\ge 1$.
	Thus, by the ``further" part of Proposition~\ref{prop:clique size}, the part (c) is true.
\end{proof}

\begin{figure}
	\begin{center}
		\subfigure[\empty]{
	\resizebox{0.25\textwidth}{!}{%
	\begin{tikzpicture}[scale=1]
	\tikzset{mynode/.style={inner sep=2pt,fill,outer sep=2.3pt,circle}}
	\node [mynode] (u1) at (-2,0) [label=below:] {};
	\node [mynode] (u2) at (0,0) [label=below :] {};
	\node [mynode] (u3) at (2,0) [label=below :] {};
	\node [mynode] (v1) at (-2,2) [label=above :] {};
	\node [mynode] (v2) at (-0,2) [label=above :] {};
	\node [mynode] (v3) at (2,2) [label=above :] {};
	\draw[->, thick] (u1) edge (v1);
	\draw[->, thick] (u1) edge (v3);
	\draw[->, thick] (u2) edge (v1);
	\draw[->, thick] (u2) edge (v2);
	\draw[->, thick] (u3) edge (v2);	
	\draw[->, thick] (u3) edge (v3);
	\end{tikzpicture}
	}
	}
	\hspace{1cm}
		\subfigure[\empty]{
	\resizebox{0.25\textwidth}{!}{%
	\begin{tikzpicture}[scale=1]
	\tikzset{mynode/.style={inner sep=2pt,fill,outer sep=2.3pt,circle}}
	\node [mynode] (u1) at (-2,0) [label=below:] {};
	\node [mynode] (u2) at (0,0) [label=below :] {};
	\node [mynode] (u3) at (2,0) [label=below :] {};
	\node [mynode] (v1) at (0,2) [label=above :] {};
	\draw[->, thick] (u1) edge (v1);
	\draw[->, thick] (u2) edge (v1);
	\draw[->, thick] (u3) edge (v1);
	\end{tikzpicture}
	}
	}	
	\hspace{1cm}
		\subfigure[\empty]{
	\resizebox{0.25\textwidth}{!}{%
	\begin{tikzpicture}[scale=1]
	\tikzset{mynode/.style={inner sep=2pt,fill,outer sep=2.3pt,circle}}
	\node [mynode] (u1) at (-2,0) [label=below:] {};
	\node [mynode] (u2) at (0,0) [label=below :] {};
	\node [mynode] (u3) at (2,0) [label=below :] {};
	\node [mynode] (v1) at (-1,2) [label=above :] {};
	\node [mynode] (v3) at (1,2) [label=above :] {};
	\draw[->, thick] (u1) edge (v1);
	\draw[->, thick] (u1) edge (u2);
	\draw[->, thick] (u2) edge (v1);
	\draw[->, thick] (u2) edge (v3);
	\draw[->, thick] (u3) edge (u2);	
	\draw[->, thick] (u3) edge (v3);
	\end{tikzpicture}
	}
	}
	\hspace{1cm}
		\subfigure[\empty]{
	\resizebox{0.14\textwidth}{!}{%
	\begin{tikzpicture}[scale=1]
	\tikzset{mynode/.style={inner sep=2pt,fill,outer sep=2.3pt,circle}}
	\node [mynode] (u1) at (-1,0) [label=below:] {};
	\node [mynode] (u2) at (0,2) [label=below :] {};
	\node [mynode] (u3) at (1,0) [label=below :] {};
	\draw[<->, thick] (u1) edge (u2);
	\draw[<->, thick] (u1) edge (u3);
	\draw[<->, thick] (u2) edge (u3);
	\end{tikzpicture}
	}
	}	
	\hspace{3cm}
		\subfigure[\empty]{
	\resizebox{0.14\textwidth}{!}{%
	\begin{tikzpicture}[scale=1]
	\tikzset{mynode/.style={inner sep=2pt,fill,outer sep=2.3pt,circle}}
	\node [mynode] (u1) at (-1,0) [label=below:] {};
	\node [mynode] (v1) at (-1,2) [label=below :] {};	
	\node [mynode] (v2) at (1,2) [label=below :] {};
	\node [mynode] (u2) at (1,0) [label=below :] {};
	\draw[<->, thick] (u1) edge (u2);
	\draw[->, thick] (u2) edge (v1);
	\draw[->, thick] (u1) edge (v1);
	\draw[<-, thick] (u2) edge (v2);
	\draw[<-, thick] (u1) edge (v2);
	\end{tikzpicture}
	}
	}	
	\end{center}
	\caption{Triangularize}
	\label{fig:triangularize}
	\end{figure}
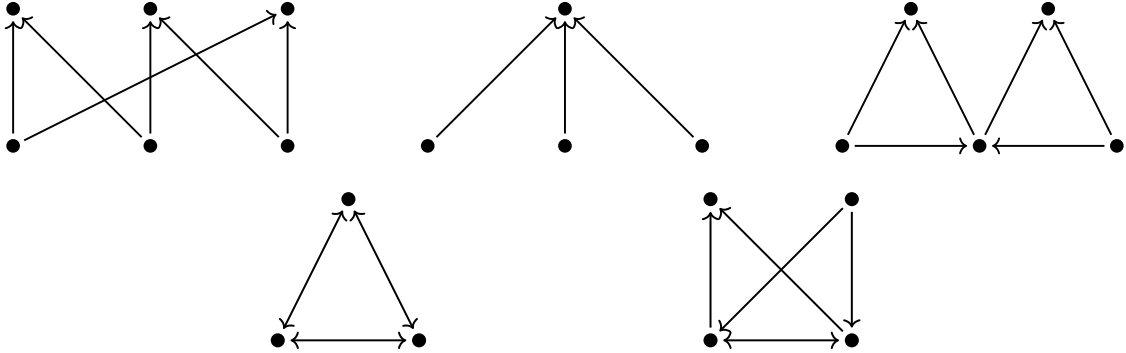

One may check that for the competition graph $G$ of a digraph $D$, $G$ has a triangle if and only if $D$ has a subdigraph isomorphic to one of the five digraphs in Figure~\ref{fig:triangularize}.
Thus we say that a digraph $D$ {\it induces a triangle} if $D$ has a subdigraph isomorphic to one of the five digraphs in Figure~\ref{fig:triangularize}.

Hefner {\it et al.} introduced the notion of irredundant digraphs.
A digraph $D$ is {\it irredundant} if $D$ has no subdigraph isomorphic to $P(2,2)$ (see Figure~\ref{fig:P(2,2)}). 
Note that any two vertices of an irredundant digraph have at most one common out-neighbor.

We say that a digraph is an \textit{$\langle \bar{i},\bar{j} \rangle$ digraph} if it is a loopless digraph satisfying $d^{-}(x)=0$ or $i$, and $d^{+}(x)=0$ or $j$ for every vertex $x$.
According to this definition, an $\langle \bar{i},\bar{j} \rangle$ digraph is an $\langle i,j \rangle$ digraph.

\begin{figure}
	\begin{center}
	\resizebox{0.15\textwidth}{!}{%
	\begin{tikzpicture}[scale=1]
	\tikzset{mynode/.style={inner sep=2pt,fill,outer sep=2.3pt,circle}}
	\node [mynode] (u) at (-1,0) [label=above:] {};
	\node [mynode] (v) at (1,0) [label=above :] {};
	\node [mynode] (x) at (-1,2) [label=above :] {};
	\node [mynode] (y) at (1,2) [label=below :] {};
	\draw[->, thick] (u) edge (x);
	\draw[->, thick] (u) edge (y);
	\draw[->, thick] (v) edge (x);
	\draw[->, thick] (v) edge (y);
	\end{tikzpicture}
	}
	\end{center}
	\caption{$P(2,2)$}
	\label{fig:P(2,2)}
	\end{figure}

\begin{proposition}\label{prop:chordal if}
	Let $G$ be the competition graph of an $\langle i,j \rangle$ digraph $D$ for some positive integers $i$ and $j$.
	Then $G$ is chordal if $D$ has no $\langle \bar{2}, \bar{2} \rangle$ subdigraph $D'$ satisfying the following: 
	\begin{enumerate}[(i)]
		\item $D'$ is irredundant;
		\item $D'$ does not induce a triangle;
		\item $D'$ has an arc.
	\end{enumerate}
\end{proposition}
\begin{proof}
	To show the contrapositive,
	suppose that $G$ is not chordal.
	That is, $G$ has a hole $H:=v_1v_2\cdots v_mv_1$ for some integer $m\ge 4$.
	Then there are vertices $u_1, u_2, \ldots, u_{m}$ in $D$ such that $u_i$ is a common out-neighbor of $v_i$ and $v_{i+1}$ for each $1\le i \le m$, where we identify $v_{m+1}$ with $v_1$.
	Now, we consider a subdigraph $D'$ of $D$ with 
	\[V(D')=\{v_1, v_2, \ldots, v_m\} \cup \{u_1, u_2, \ldots, u_{m}\}\]
	and 
	\[A(D')=\{(v_i,u_i), (v_{i+1},u_i) \colon\, 1\le i \le m\}.\]
	Clearly, $D'$ has an arc and it is irredundant. 
	If $u_i=u_j$ for some distinct $i$ and $j$, then $v_i$, $v_{i+1}$, $v_j$, and $v_{j+1}$, at least three of which are distinct, form a clique in $G$, which contradicts to the fact that $H$ is a hole. 
	Thus $u_1, u_2, \ldots, u_{m}$ are all distinct.
	Therefore $D'$ is a $\langle \bar{2}, \bar{2}\rangle$ digraph and the competition graph of $D'$ is $H$ with isolated vertices.
	If $D'$ induces a triangle, then the competition graph of $D'$ contains a triangle, which is impossible.
	Thus $D'$ does not induce a triangle.
\end{proof}

Especially, the sufficient condition given in Proposition~\ref{prop:chordal if} becomes a necessary and sufficient condition when $j=2$ (Theorem~\ref{thm:i,2 chordal}).

To avoid unnecessary repetition, we call a $\langle \bar{2}, \bar{2} \rangle$ digraph satisfying the conditions (i), (ii), and (iii)  in Proposition~\ref{prop:chordal if}  a $\langle \bar{2}, \bar{2} \rangle$ {\it  good} digraph.

\begin{theorem}\label{thm:i,2 chordal}
	Let $G$ be the competition graph of an $\langle i,2 \rangle$ digraph $D$ for some positive integer $i$.
	Then $G$ is chordal if and only if $D$ has no $\langle \bar{2}, \bar{2} \rangle$ good subdigraph.
	\end{theorem}
\begin{proof}
	By Proposition~\ref{prop:chordal if}, it suffices to show the ``only if" part.
	To the contrary,
	suppose that $G$ is chordal while $D$ has a $\langle \bar{2}, \bar{2} \rangle$ good subdigraph $D'$. 
	Let $G'$ be the competition graph of $D'$.
	To show that a vertex has degree $0$ or $2$ in $G'$, take a vertex $v$ in $D'$.
	If $d^+_{D'}(v)=0$, then $v$ is an isolated vertex in $G'$.
	Suppose $d^+_{D'}(v)\neq 0$.
	Since $D'$ is a $\langle \bar{2}, \bar{2} \rangle$ digraph, $d^+_{D'}(v)=2$ and so there are out-neighbors $x$ and $y$ of $v$.
	For the same reason, $x$ (resp.\ $y$) has an in-neighbor $u$ (resp.\ $w$) which is distinct from $v$.
	Then $u$ and $w$ are neighbors of $v$ in $G$.
	Since $D'$ is irredundant, $u$ and $w$ are distinct.
	Thus $v$ has at least two neighbors in $G'$.
	By Corollary~\ref{cor:edge clique cover}(c), $v$ has degree $2$ in $G'$.
	Since $v$ is arbitrarily chosen, we have shown that each vertex  has degree $0$ or $2$ in $G'$.
	Since $D'$ has an arc, it has a vertex of indegree $2$ and so $G'$ has an edge.
	Thus $G'$ has a cycle.
	Since $D'$ does not induce a  triangle, no cycle of $G'$ is a triangle.
	Therefore $G'$ has a hole $H:=v_1v_2\cdots v_mv_1$ for some integer $m\ge 4$.
	Then there are vertices $u_1, u_2, \ldots, u_{m}$ in $D'$ such that $u_i$ is a common out-neighbor of $v_i$ and $v_{i+1}$ for each $1\le i \le m$ where we identify $v_{m+1}$ with $v_1$.
	Since $D'$ is a $\langle \bar{2}, \bar{2} \rangle$ digraph, $u_1, u_2, \ldots, u_{m}$ are all distinct.
	Thus, for each $1\le i \le m$,
	\[N^+_{D'}(v_i)=\{u_{i-1}, u_i\}\]
	where we identify $u_0$ with $u_m$.
	Then, since $D$ is an $\langle i,2 \rangle$ digraph,
	\[N^+_{D}(v_i)=N^+_{D'}(v_i)=\{u_{i-1}, u_i\}  \]
	for each $1\le i \le m$.
	We note that $\{u_{i-1},u_i\} \cap \{u_{j-1}, u_j\}=\emptyset$ if $|i-j| \ge 2$.
	Therefore, for any nonconsecutive $v_i$ and $v_j$ in $H$, they have no common out-neighbor in $D$ and so they are not adjacent in $G$.
	Hence, $H$ is a hole of $G$, which is a contradiction that $G$ is chordal.	
	\end{proof}

\begin{corollary}
	Let $G$ be the competition graph of a $\langle 2,2 \rangle$ digraph $D$.
	Then $G$ is interval if and only if there is no $\langle \bar{2}, \bar{2} \rangle$ good subdigraph of $D$.
\end{corollary}
\begin{proof}
	By Corollary~\ref{cor:edge clique cover}(c), each vertex of $G$ has degree at most $2$.
	Thus $G$ is the union of paths and cycles.
	Therefore $G$ is interval if and only if it is chordal, and so the statement is immediately true by Theorem~\ref{thm:i,2 chordal}.
\end{proof}

\section{The set containment relation on $\{\mathcal{G}_{\langle i,j \rangle}\colon\, i,j \ge 1\}$}\label{sec:3}

We recall that $\mathcal{G}_{\langle i,j \rangle}$ denotes the family of $\langle i,j \rangle$ competition graphs.
In this section, we take a look at set containment between $\mathcal{G}_{\langle i,j \rangle}$ and $\mathcal{G}_{\langle k,l \rangle}$ for given positive integers $i$, $j$, $k$, and $l$ with $(i,j) \neq (k,l)$.
By symmetry, it is sufficient to consider the following four cases.
\begin{center}
	(1) $i=k$, $j>l$; \quad (2) $i>k$, $j>l$; \quad (3) $i>k$, $j=l$; \quad (4) $i>k$, $j<l$.
\end{center}
This section will focus on proving the following theorem.

\begin{theorem}\label{thm:containment}
	Let $i$, $j$, $k$, and $l$ be positive integers satisfying $(i,j) \neq (k,l)$.
	Without loss of generality, assume $i\ge k$.
	Then there are four cases and the following hold.
	\begin{enumerate}[(1)]
		\item If $i=k$ and $j>l$, then \[\begin{cases}
			\mathcal{G}_{\langle i,j \rangle}=\mathcal{G}_{\langle k,l \rangle} & \text{if }i=k=1;\\
			\mathcal{G}_{\langle k,l \rangle} \subsetneq \mathcal{G}_{\langle i,j \rangle} & \text{if }i=k>1.
		\end{cases} \]
		\item If $i>k$ and $j>l$, then $\mathcal{G}_{\langle k,l \rangle} \subsetneq \mathcal{G}_{\langle i,j \rangle}$.
		\item If $i>k$ and $j=l$, then $\mathcal{G}_{\langle k,l \rangle} \subsetneq \mathcal{G}_{\langle i,j \rangle}$.
		\item If $i>k$ and $j<l$, then the following are true.
		\begin{enumerate}[(a)]
		\item \[\begin{cases}
			\mathcal{G}_{\langle k,l \rangle} \subsetneq \mathcal{G}_{\langle i,j \rangle} & \text{if }k=1;\\
			\mathcal{G}_{\langle k,l \rangle} \not\subseteq \mathcal{G}_{\langle i,j \rangle} & \text{if }k>1.
		\end{cases}\]
		\item $\mathcal{G}_{\langle i,j \rangle}\not\subseteq \mathcal{G}_{\langle k,l \rangle}$ if $i>k^2-k+1$ or $(i-1)j>(k-1)l$ or $j\ge k\ge 2$.
		\end{enumerate} 
		\end{enumerate}
\end{theorem}

In order to establish the proof, we undertake the following necessary preparations.

\begin{proposition}\label{prop:1,j}
	A graph $G$ is a $\langle 1,j \rangle$ competition graph for some positive integer $j$ if and only if $G$ is edgeless, that is, $E(G)=\emptyset$.
\end{proposition}
\begin{proof}
An edgeless graph is an $\langle i,j \rangle$ digraph for any positive integers $i$ and $j$ and is the competition graph of itself.
Thus the ``if" part is clear.
The ``only if" part directly follows Corollary~\ref{cor:edge clique cover}(a).
\end{proof}
	
\begin{proposition}\label{prop:i,1}
	A graph $G$ is an $\langle i,1 \rangle$ competition graph for some positive integer $i$ if and only if $G$ is trivial or is the disjoint union of at least two complete graphs each of which has size at most $i$. 
\end{proposition}
\begin{proof}
	To show the ``only if" part, suppose that $G$ is an $\langle i,1 \rangle$ competition graph for some positive integer $i$.
	By Theorem~\ref{thm:edge clique cover_2}(iii), there is an edge clique cover $\mathcal{C}$ in $\mathbf{C}(G,i,1)$.
	Then each vertex of $G$ belongs to at most one clique in $\mathcal{C}$ and so $G$ is the disjoint union of cliques each of which has size at most $i$.
	Moreover, by Theorem~\ref{thm:edge clique cover_2}(ii), $G$ is not a nontrivial complete graph.
	Thus $G$ is trivial or is the disjoint union of at least two complete graphs each of which has size at most $i$. 
	
	To show the ``if" part, suppose that $G$ is trivial or is the disjoint union of at least two complete graphs each of which has size at most $i$ for some positive integer $i$. 
	Then $G$ is not a nontrivial complete graph and the collection of components of $G$ is an edge clique cover in $\mathbf{C}(G,i,1)$.
	Thus $G$ is an $\langle i,1 \rangle$ competition graph by Theorem~\ref{thm:edge clique cover_2}.
\end{proof}

\begin{proposition}\label{prop:tool1}
	Let $i$, $j$, $k$, $l$, $m$, and $n$ be positive integers satisfying $i,j \ge 2$, $j>l$, and $(i-1)j>(m-1)n$. 
	Then $\mathcal{G}_{\langle i,j \rangle}\not\subseteq \mathcal{G}_{\langle k,l \rangle}\cup \mathcal{G}_{\langle m,n \rangle}$.
\end{proposition}
\begin{proof}
	Consider a star $K_{1,j}$ with the center $u$ and the leaves $u_1, \ldots, u_j$.
	We replace each edge $uu_t$ of $K_{1,j}$ with the complete graph $C_t$ having $\{u, u^1_t, \ldots, u^{i-1}_t\}$ as the vertex set  for each $1\le t \le j$ to obtain a graph $G^*_{i,j}$ (see Figure~\ref{fig:G^*_{5,3}} for $G^*_{5,3}$).
	Then $\{C_1, \ldots, C_j\} \in \mathbf{C}(G^*_{i,j},i,j)$.
	Since $i,j \ge 2$, $G^*_{i,j}\not\cong K_2$.
	Thus $G^*_{i,j} \in \mathcal{G}_{\langle i,j \rangle}$ by Theorem~\ref{thm:edge clique cover_2}.
	We note that the subgraph of $G^*_{i,j}$ induced by $\{u,u_1^1, \ldots, u_j^1\}$ is isomorphic to $K_{1,j}$.
	Then, since $j>l$, $G^*_{i,j}$ has a subgraph isomorphic to $K_{1,l+1}$ and so, by Corollary~\ref{cor:edge clique cover}(b), 
	$G^*_{i,j} \not\in \mathcal{G}_{\langle k,l \rangle}$.
	Since $u$ has degree $(i-1)j$ in $G^*_{i,j}$ and $(i-1)j>(m-1)n$, $G^*_{i,j} \not\in \mathcal{G}_{\langle m,n \rangle}$ by Corollary~\ref{cor:edge clique cover}(c).
	Therefore $G^*_{i,j} \in \mathcal{G}_{\langle i,j \rangle}
	- (\mathcal{G}_{\langle k,l \rangle}\cup \mathcal{G}_{\langle m,n \rangle})$ and so $\mathcal{G}_{\langle i,j \rangle}\not\subseteq \mathcal{G}_{\langle k,l \rangle}\cup \mathcal{G}_{\langle m,n \rangle}$.
\end{proof}

\begin{figure}
	\begin{center}
	\resizebox{0.4\textwidth}{!}{%
	\begin{tikzpicture}[scale=1]
	\tikzset{mynode/.style={inner sep=2pt,fill,outer sep=2.3pt,circle}}
	\node [mynode] (u) at (0,0) [label=above:$u$] {};
	\node [mynode] (u11) at (-25:2cm) [label=right: $u_1^1$] {};
	\node [mynode] (u12) at (5:2.5cm) [label=right: $u_1^2$] {};
	\node [mynode] (u13) at (35:2.5cm) [label=right: $u_1^3$] {};
	\node [mynode] (u14) at (65:2cm) [label=above: $u_1^4$] {};
	\node [mynode] (u21) at (115:2cm) [label=above: $u_2^1$] {};
	\node [mynode] (u22) at (145:2.5cm) [label=above: $u_2^2$] {};
	\node [mynode] (u23) at (175:2.5cm) [label=left: $u_2^3$] {};
	\node [mynode] (u24) at (205:2cm) [label=left: $u_2^4$] {};
	
	\node [mynode] (u31) at (225:2cm) [label=left: $u_3^1$] {};
	\node [mynode] (u32) at (255:2.5cm) [label=below: $u_3^2$] {};
	\node [mynode] (u33) at (285:2.5cm) [label=below: $u_3^3$] {};
	\node [mynode] (u34) at (315:2cm) [label=right: $u_3^4$] {};	
	
	\draw[-, thick] (u) edge (u11);	
	\draw[-, thick] (u) edge (u12);
	\draw[-, thick] (u) edge (u13);
	\draw[-, thick] (u) edge (u14);
	\draw[-, thick] (u11) edge (u12);	
	\draw[-, thick] (u11) edge (u13);	
	\draw[-, thick] (u11) edge (u14);
	\draw[-, thick] (u12) edge (u13);
	\draw[-, thick] (u12) edge (u14);
	\draw[-, thick] (u13) edge (u14);

	\draw[-, thick] (u) edge (u21);
	\draw[-, thick] (u) edge (u22);
	\draw[-, thick] (u) edge (u23);
	\draw[-, thick] (u) edge (u24);
	\draw[-, thick] (u21) edge (u22);	
	\draw[-, thick] (u21) edge (u23);	
	\draw[-, thick] (u21) edge (u24);
	\draw[-, thick] (u22) edge (u23);
	\draw[-, thick] (u22) edge (u24);
	\draw[-, thick] (u23) edge (u24);
	
	\draw[-, thick] (u) edge (u31);	
	\draw[-, thick] (u) edge (u32);
	\draw[-, thick] (u) edge (u33);
	\draw[-, thick] (u) edge (u34);
	\draw[-, thick] (u31) edge (u32);	
	\draw[-, thick] (u31) edge (u33);	
	\draw[-, thick] (u31) edge (u34);
	\draw[-, thick] (u32) edge (u33);
	\draw[-, thick] (u32) edge (u34);
	\draw[-, thick] (u33) edge (u34);
	\end{tikzpicture}
	}
	\end{center}
	\caption{$G^*_{5,3}$}
	\label{fig:G^*_{5,3}}
	\end{figure}

\begin{proposition}\label{prop:tool2}
	Let $i$, $j$, $k$, and $l$ be positive integers satisfying $i>k^2-k+1$. 
	Then $\mathcal{G}_{\langle i,j \rangle}\not\subseteq \mathcal{G}_{\langle k,l \rangle}$.
\end{proposition}
\begin{proof}
	By Proposition~\ref{prop:i,1}, $K_i \cup K_i \in \mathcal{G}_{\langle i,1 \rangle} \subseteq \mathcal{G}_{\langle i,j \rangle}$. 
	If $K_i \cup K_i \in \mathcal{G}_{\langle k,l \rangle}$, then, by Corollary~\ref{cor:edge clique cover}(a), 
	\[ (i-1)i=|E(K_i \cup K_i)| \le \frac{k(k-1)}{2}|V(K_i \cup K_i)|=(k^2-k)i \]
	and so $i\le k^2-k+1$.
	Since $i>k^2-k+1$, $K_i \cup K_i \not\in \mathcal{G}_{\langle k,l \rangle}$.
	Thus $K_i \cup K_i \in \mathcal{G}_{\langle i,j \rangle}-\mathcal{G}_{\langle k,l \rangle}$.
	\end{proof}

\begin{proposition}\label{prop:tool3}
	Let $i$, $j$, $k$, and $l$ be positive integers satisfying $j\ge k\ge 2$ and $i>k$. 
	Then $\mathcal{G}_{\langle i,j \rangle}\not\subseteq \mathcal{G}_{\langle k,l \rangle}$.
\end{proposition}
\begin{proof}
	It suffices to construct a $\langle k+1,k \rangle$ competition graph $G$ satisfying 
	\[\frac{|E(G)|}{|V(G)|} > \frac{k(k-1)}{2}. \]
	For, $G\in\mathcal{G}_{\langle k+1,k \rangle} \subseteq \mathcal{G}_{\langle i,j \rangle}$ while $G \notin \mathcal{G}_{\langle k,l \rangle}$ by Corollary~\ref{cor:edge clique cover}(a).

	Given positive integers $m$ and $n$, we consider the Cartesian product of $n$ copies of $[m]=\{1, \ldots, m\}$, i.e. 
	\[[m]^n:= \{ (x_1, \ldots, x_n) \in \mathbb{Z}^n \colon\, 1\le x_1, \ldots, x_n \le m \}. \]
	Now we consider the graph $G$ such that $V(G)=[k+1]^k$
	and two vertices are adjacent if and only if they differ in exactly one component.
	Then $|V(G)|=(k+1)^k$.
	Moreover, each vertex is adjacent with exactly other $k^2$ vertices and so \[ |E(G)|=\frac{1}{2}\sum_{v\in V(G)}d(v)=\frac{k^2|V(G)|}{2}.\]
	Thus $G$ satisfies 
	\[\frac{|E(G)|}{|V(G)|} > \frac{k(k-1)}{2}. \]
	It remains to show  $G \in \mathcal{G}_{\langle k+1,k \rangle}$.
	Given an integer  $t \in [k]$ and a vertex $u$ of $G$, 
	let 
	\[C_t^u=\{(x_1, \ldots, x_k) \in V(G) \colon\, x_i=u_i \text{ for each } i\in [k]-\{t\}\}. \]
	Then each $C_t^u$ forms a clique of size $k+1$ in $G$.
	If $u$ and $v$ are adjacent in $G$, then they differ in exactly one component, say the $p$th component, and so they are contained in $C_p^u$.
	Thus $\mathcal{C}:=\{C_t^u \colon\,1\le t \le k, u\in [k+1]^{k}\}$ is an edge clique cover of $G$.
	For each vertex $u$ of $G$, there are exactly $k$ cliques $C_1^u, \ldots, C_k^u$ in $\mathcal{C}$ containing $u$.
	Thus 
	\[ (k+1)|\mathcal{C}|=\sum_{C \in \mathcal{C}}\sum_{u\in C} 1=|\{(C,u) \colon\, C \in \mathcal{C}, u \in C \}|= \sum_{u\in V(G)}\sum_{u\in C \in \mathcal{C}} 1 =k|V(G)|\]
	and so $|\mathcal{C}| \le |V(G)|$.
	Therefore \[\mathcal{C} \in \mathbf{C}(G,k+1,k)\]
	and so, by Theorem~\ref{thm:edge clique cover_2}, $G$ is a $\langle k+1, k \rangle$ competition graph.
\end{proof}

Now we are ready to prove Theorem~\ref{thm:containment}.

\begin{proof}[Proof of Theorem~\ref{thm:containment}]
	To show (1), suppose that $i=k$ and $j>l$.
	Then $\mathcal{G}_{\langle k,l \rangle} \subseteq \mathcal{G}_{\langle i,j \rangle}$ and $j \ge 2$.
	If $i=k=1$, then $\mathcal{G}_{\langle i,j \rangle}=\mathcal{G}_{\langle k,l \rangle}$ by Proposition~\ref{prop:1,j}.
	If $i=k > 1$, then $\mathcal{G}_{\langle i,j \rangle} \not\subseteq \mathcal{G}_{\langle k,l \rangle}$ by Proposition~\ref{prop:tool1} and so $\mathcal{G}_{\langle k,l \rangle} \subsetneq \mathcal{G}_{\langle i,j \rangle}$.
	Thus (1) holds.

	Suppose that $i>k$ and $j>l$.
	Then $\mathcal{G}_{\langle k,l \rangle} \subseteq \mathcal{G}_{\langle i,j \rangle}$ and $i,j \ge 2$.
	Moreover, $\mathcal{G}_{\langle k,l \rangle} \subsetneq \mathcal{G}_{\langle i,j \rangle}$ by Proposition~\ref{prop:tool1}.
	Thus (2) holds.
	
	Suppose that $i>k$ and $j=l$.
	Then $\mathcal{G}_{\langle k,l \rangle} \subseteq \mathcal{G}_{\langle i,j \rangle}$ and $i \ge 2$.
	If $j=l=1$, then the disjoint union of two complete graphs each of whose order is $i$ is an $\langle i,j \rangle$ competition graph which is not a $\langle k,l \rangle$ competition graph by Proposition~\ref{prop:i,1} and so $\mathcal{G}_{\langle k,l \rangle} \subsetneq \mathcal{G}_{\langle i,j \rangle}$.
	If $j=l\ge2$, then Proposition~\ref{prop:tool1} implies $\mathcal{G}_{\langle i,j \rangle} \not\subseteq \mathcal{G}_{\langle k,l \rangle}$ by substituting $(i,j,k,l)$ for $(i,j,m,n)$.
	Thus $\mathcal{G}_{\langle k,l \rangle} \subsetneq \mathcal{G}_{\langle i,j \rangle}$.
	Therefore (3) holds.	
	
	Finally, to show (4), suppose that $i>k$ and $j<l$.
	Then $i \ge 2$ and so it is easy to check that $K_2 \cup K_1 \in \mathcal{G}_{\langle i,j \rangle}$ by Theorem~\ref{thm:edge clique cover_2}.
	Thus, if $k=1$, $\mathcal{G}_{\langle k,l \rangle} \subsetneq \mathcal{G}_{\langle i,j \rangle}$ by Propositions~\ref{prop:1,j}.
	If $k \ge 2$, then Proposition~\ref{prop:tool1} implies $\mathcal{G}_{\langle k,l \rangle} \not\subseteq \mathcal{G}_{\langle i,j \rangle}$ by substituting $(k,l,i,j)$ for $(i,j,k,l)$.
	Thus (a) is true.
	Now we prove (b).
	If $i>k^2-k+1$, then $\mathcal{G}_{\langle i,j \rangle}\not\subseteq \mathcal{G}_{\langle k,l \rangle}$ by Proposition~\ref{prop:tool2}.
	Suppose $(i-1)j>(k-1)l$. If $k=1$, then we have already checked that $\mathcal{G}_{\langle k,l \rangle} \subsetneq \mathcal{G}_{\langle i,j \rangle}$ and so $\mathcal{G}_{\langle i,j \rangle}\not\subseteq \mathcal{G}_{\langle k,l \rangle}$.
	If $k\ge 2$, then $\mathcal{G}_{\langle i,j \rangle}\not\subseteq \mathcal{G}_{\langle k,l \rangle}$ by Proposition~\ref{prop:tool1}. 
	If $j\ge k\ge 2$, then $\mathcal{G}_{\langle i,j \rangle}\not\subseteq \mathcal{G}_{\langle k,l \rangle}$ by Proposition~\ref{prop:tool3}.
	Hence we have shown that (4) holds.
\end{proof}

\section{Concluding remarks}
To fully address the case (4) of Theorem~\ref{thm:containment}, we need to find out whether or not there is a graph belonging to $\mathcal{G}_{\langle i,j \rangle}- \mathcal{G}_{\langle k,l \rangle}$ when $k^2-k+1 \ge i>k \ge 2$, $j<\min(k,l)$, and $(i-1)j\le(k-1)l$.
This problem seems challenging to resolve due to the following reasons.
Consider $(i,j,k,l)=(3,1,2,t)$ and $(i,j,k,l)=(4,1,3,t)$ for each integer $t\ge 2$, and $(i,j,k,l)=(5,1,3,2)$ all of which satisfy the conditions given above.
Yet, one may check that $\mathcal{G}_{\langle 3,1 \rangle}\subsetneq \mathcal{G}_{\langle 2,l \rangle}$ and $\mathcal{G}_{\langle 4,1 \rangle}\subsetneq \mathcal{G}_{\langle 3,l \rangle}$ for all $l\ge 2$ by Proposition~\ref{prop:i,1} and Theorem~\ref{thm:edge clique cover_2}, while
$ \mathcal{G}_{\langle 5,1 \rangle} \not\subseteq \mathcal{G}_{\langle 3,2 \rangle}$ since
$K_5 \cup K_5 \in \mathcal{G}_{\langle 5,1 \rangle}- \mathcal{G}_{\langle 3,2 \rangle}$.

\section{Acknowledgement}
This work was supported by Science Research Center Program through the National Research Foundation of Korea(NRF) Grant funded by the Korean Government (MSIP)(NRF-2022R1A2C\\1009648 and 2016R1A5A1008055).
Especially, the first author was supported by Basic Science Research Program through the National Research Foundation of Korea(NRF) funded by the Ministry of Education (NRF-2022R1A6A3A13063000).


\end{document}